\documentclass[12pt]{amsart}

\usepackage{amsmath}
\usepackage{amsfonts}
\usepackage{amssymb}

 \newtheorem{theorem}{Theorem}

\author{Sergey V. Galaev}

\title{The Intrinsic Geometry of Almost Contact Metric Manifolds}

\begin{document}

\maketitle

\begin{abstract}
 In this paper the notion  of the intrinsic
geometry of an almost contact metric manifold is introduced.
Description of some classes of spaces with almost contact metric
structures in terms of the intrinsic geometry is given. A new type
of almost contact metric spaces, more precisely, Hermitian almost
contact metric spaces, is introduced. \end{abstract}

\textbf{Key words:} almost contact manifold, Sasakian manifold, $K$-contact manifold, the intrinsic geometry of almost contact metric manifolds.



\section*{Introduction}

The research of the geometry of manifolds with almost contact
metric structures has been begun in the fundamental papers by
Chern \cite{Lit1}, J.\,Gray \cite{Lit2} and Sasaki \cite{Lit3}.
Almost contact metric structures constitute the odd-dimensional
analog of almost Hermitian structures. There are a lot of
important interplays between these structures. In the same time,
the geometry of almost contact metric structures is appreciably
different from the geometry of almost Hermitian structures and its
study requires in principle new tools. The results obtained in
this area before 1976 are reflected  in full measure in the book
\cite{Lit4}. Important contribution to the development of the
geometry of almost contact metric spaces inserted V.F. Kirichenko
and his students, see \cite{Lit5,Lit6}. One can get the impression
of the last achievements in this area and about applications to
theoretical physics from the works \cite{Lit7,Lit8}.

In the present paper, the notion of the intrinsic  geometry of a
manifold with an almost contact metric structure is introduced. In
the terminology developed by V.V.Wagner \cite{Lit10}, the manifold
with an almost contact metric structure is a nonholonomic manifold
of codimension 1 with additional structures. These structures
Wagner called intrinsic. The notion of the intrinsic geometry of a
nonholonomic manifold was defined by Schouten  as the properties
that depend only on the parallel transport in the nonholonomic
manifold and on the closing of the nonholonomic manifold in the
ambient manifold.

In \cite{Lit9} Wagner writes: "Schouten shows the possibility of
the direct construction of the intrinsic geometry of a
nonholonomic manifold without usage of the parallel transport in
the ambient space. In the case of metric nonholonomic manifold,
the intrinsic geometry is defined by assigning a closing and a
metric on local tangent spaces". Developing the intrinsic geometry
of a nonholonomic manifold, Wagner defines and investigates the
curvature tensor of a nonholonomic manifold that generalizes the
Schouten curvature tensor. The curvature tensor, which later was
called the Wagner curvature tensor, first was constructed for a
nonholonomic metric manifold of arbitrary codimension \cite{Lit9},
and then this result was specified for the case of a nonholonomic
manifold of codimension 1 endowed with an intrinsic linear
connection \cite{Lit10}. Wagner used his theory of the curvature
of  nonholonomic manifolds for solving some problems of classical
mechanics and calculus of variations. In the present paper we
propose to use the methods of nonholonomic geometry developed by
Wagner for investigation the geometry of manifolds with almost
contact metric structure. The new approach allows to pick out new
types of spaces. For example, we give the definition of an
Hermitian almost contact metric spaces. The already known results
obtain new description on the language of the intrinsic geometry.
Following the ideology developed in the works of Schouten and
Wagner, we define the intrinsic geometry of an almost contact
metric space $X$ as the aggregate of the properties that possess
the following objects: a smooth distribution $D$ defined by a
contact form   $\eta$; an admissible field of endomorphisms
$\varphi$ of $D$ (which we call an admissible almost complex
structure) satisfying $\varphi^{2}=-1$; an admissible Riemannian
metric field $g$ that is related to $\varphi$ by
$g(\varphi\vec{X},\varphi\vec{Y})=g(\vec{X},\vec{Y})$, where
$\vec{X}$ and $\vec{Y}$ are admissible vector fields. To the
objects of the intrinsic geometry of an almost contact metric
space one should ascribe also the  objects derived from the just
mentioned: the 2-form $\omega=d\eta$; the vector field $\vec{\xi}$
(which is called the Reeb vector field) defining the closing
$D^{\bot}$ of $D$, i.e. $\vec{\xi} \in D^{\bot}$, and given by the
equalities $\eta(\vec{\xi})=1$, ${\rm ker\,} \omega={\rm span}
(\vec{\xi})$ in the case when the 2-form  $\omega$ is of maximal
rank; the intrinsic connection $\nabla$ that  defines the parallel
transport of admissible vectors along admissible curves and is
defined by the metric $g$; the connection $\nabla^{1}$ that is a
natural extension of the connection $\nabla$ which accomplishes
the parallel transport of admissible vectors along arbitrary
curves of the manifold $X$.

The paper consists of two sections. In the first section we
provide the basic concepts of the theory of manifolds with almost
contact metric structure. We introduce the notion of the adapted
coordinate system. The adapted coordinates play in the geometry of
the nonholonomic manifolds the same role as the holonomic
coordinates on a holonomic manifold, see e.g. \cite{Lit10}. The
adapted coordinates are extensively used in the geometry of
foliations \cite{Lit11}. Next we introduce the notion of an
admissible (to the distribution $D$) tensor structure. An
admissible tensor structure is an object of the intrinsic geometry
of a nonholonomic manifold \cite{Lit10}. In the literature on the
geometry of the fibering spaces, the admissible tensor structures
are usually called semi basic. We give some information about the
intrinsic connections compatible with admissible tensor
structures. Among the connection compatible with the admissible
Riemannian metric, we study the connections compatible with an
admissible almost complex structure. We discuss the connection
over a distribution that was introduced in \cite{Lit12,Lit13} and
applied in \cite{Lit14,Lit15} to  manifolds with an admissible
Finsler metric.

In the second section we expound some of the main theses of the
geometry of almost contact metric spaces in terms of the intrinsic
geometry. It is shown that the almost contact metric structure
defined in the intrinsic way corresponds to a certain almost
contact metric structure defined in the usual way. The intrinsic
connection is used for description and characterization of the
normal and Sasakian structures.

\section{Admissible tensor structures and intrinsic connection compatible with them}

Let $X$  be a smooth manifold of  an odd dimension $n$. Denote by
$\Xi(X)$ the $C^{\infty}(X)$-module of smooth vector fields on $X$
and by $d$ the exterior derivative. All manifolds, tensors and
other geometric objects will be assumed to be smooth of the class
$C^{\infty}$. For simplification, in what follows we call tensor
fields simply be tensors. An almost contact metric structure on
 $X$ is an aggregate
 $(\varphi, \vec{\xi}, \eta, g)$ of the tensor fields on $X$, where $\varphi$ is a tensor field of type $(1, 1)$, which is called the structure
 endomorphism, $\vec{\xi}$ and $\eta$  are vector and covector,
 which are called the structure vector and the contact form,
 respectively, and  $g$ is a (pseudo-)Riemannian metric.
 Moreover,
$$\eta(\vec{\xi})=1,\quad \varphi(\vec{\xi})=0,\quad \eta \circ
\varphi=0,$$
$$\varphi^2\vec{X}=-\vec{X}+\eta(\vec{X})\vec{\xi},\quad
g(\varphi\vec{X},\varphi\vec{Y})=g(\vec{X},\vec{Y})-\eta(\vec{X})\eta(\vec{Y})$$
for all $\vec{X}, \vec{Y} \in \Xi(X)$. It is easy to check that
the tensor
 $\Omega(\vec{X}, \vec{Y})=g(\vec{X},
\varphi\vec{Y})$ is skew-symmetric. It is called the fundamental
tensor of the structure. A manifold with a fixed almost contact
metric structure is called an almost contact metric manifold. If
 $\Omega=d\eta$ holds, then the almost contact metric structure is called a contact metric structure.
An almost contact metric structure is called normal if
$$N_{\varphi}+2d\eta\otimes\vec{\xi}=0,$$ where $N_{\varphi}$ is the
Nijenhuis torsion defined for the tensor $\varphi$. A normal
contact metric structure is called a Sasakian structure. A
manifold with a given Sasakian structure is called a Sasakian
manifold. Let $D$ be the smooth distribution of codimension 1
defined by the form  $\eta$, and $D^\bot={\rm span}(\vec{\xi})$ be
the closing of $D$. In what follows we assume that the restriction
of the 2-form $\omega=d\eta$ to the distribution $D$  is
non-degenerate. In this case the vector $\vec{\xi}$ is uniquely
defined by the condition $$\eta(\vec{\xi})=1,\quad {\rm ker\,}
\omega={\rm span} (\vec{\xi}),$$ and  it is called the Reeb vector
field. The smooth distribution $D$ we call sometimes a
nonholonomic manifold.

For investigation of the intrinsic geometry of a nonholonomic
manifold, and  generally for the study of almost contact metric
structures, it is suitable to use  coordinate systems satisfying
certain additional conditions. We say that a coordinate map
$K(x^\alpha)$ $(\alpha,\beta,\gamma=1,...,n)$
$(a,b,c,e=1,...,n-1)$ on a manifold $X$ is adapted to the
nonholonomic manifold  $D$ if $$D^{\bot}={\rm
span}\left(\frac{\partial}{\partial x^{n}}\right)$$ holds. It is
easy to show that any two adapted coordinate map are related by a
transformation of the form  $$x^{a}=x^{a}(x^{\tilde{a}}),\quad
x^{n}=x^{n}(x^{\tilde{a}},x^{\tilde{n}}).$$ Such coordinate
systems are called by Wagner in \cite{Lit10} gradient coordinate
systems. Adapted coordinates have their applications in the
foliation theory, see e.g. \cite{Lit7}.

Let $P:TX\rightarrow D$ be the projection map defined by the
decomposition $TX=D\oplus D^{\bot}$ and let $K(x^{\alpha})$ be an
adapted coordinate map.  Vector fields
$$P(\partial_{a})=\vec{e}_{a}=\partial_{a}-\Gamma^{n}_{a}\partial_{n}$$
are linearly independent, and linearly generate the system $D$
over the domain of the definition of the coordinate map: $$D={\rm
span}(\vec{e}_{a}).$$ Thus we have on $X$ the nonholonomic field
of bases  $(\vec{e}_{a},\partial_{n})$ and the corresponding field
of cobases $$(dx^a,\theta^{n}=dx^{n}+\Gamma^{n}_{a}dx^{a}).$$ It
can be checked directly that
$$[\vec{e}_{a},\vec{e}_{b}]=M^{n}_{ab}\partial_{n},$$ where the
components  $M^{n}_{ab}$ form the so-called tensor of
nonholonomicity \cite{Lit10}. Under assumption that for all
adapted coordinate systems it holds  $\vec{\xi}=\partial_{n}$, the
following equality takes place
$$[\vec{e}_{a},\vec{e}_{b}]=2\omega_{ba}\partial_{n},$$ where
$\omega=d\eta$. In what follows we consider exceptionally adapted
coordinate systems that satisfy the condition
$\vec{\xi}=\partial_{n}$. We say also that the basis
$$\vec{e}_{a}=\partial_{a}-\Gamma^{n}_{a}\partial_{n}$$ is adapted,
as the basis defined by an adapted coordinate map. Under the
transformation of the adapted coordinate systems, the vectors of
the adapted bases transform in the following way: $
\vec{e}_{a}=\frac{\partial x^{\tilde{a}}}{\partial
x^{a}}\vec{e}_{\tilde{a}}.$

We call a tensor field defined on an almost contact metric
manifold admissible (to the distribution $D$) if it vanishes
whenever its vectorial argument belongs to the closing
 $D^\bot$ and its covectorian argument is proportional to the form  $\eta$. The coordinate form of an admissible tensor field
of type  $(p,q)$ in an adapted coordinate map looks like
$$
t=t^{a_{1},...,a_{p}}_{b_{1},...,b_{q}}\vec{e}_{a_{1}}\otimes...\otimes\vec{e}_{a_{p}}\otimes dx^{b_{1}}\otimes...\otimes dx^{b_{q}}.
$$
In particular, an admissible vector field is a vector field
tangent to the distribution $D$, and an admissible 1-form is a
1-form zero on the closing  $D^\bot$. It is clear that any tensor
structure defined on the manifold  $X$ defines on it a unique
admissible tensor structure of the same type. From the definition
of an almost contact structure it follows that the field of
endomorphisms $\varphi$ is an admissible tensor field of type $(1,
1)$. The field of endomorphisms  $\varphi$ we call an admissible
almost complex structure. The 2-form  $\omega=d\eta$ is also an
admissible tensor field. In the geometry of the fibered spaces an
admissible tensor field is called semi basic.

\begin{theorem}\label{Th1} The derivatives
  $\partial_{n}t$ of the components of an admissible tensor field $t$
  in an adapted coordinate system are components of an admissible tensor
  field of the same type. \end{theorem}

  The proof of the theorem follows from the fact that the
  components of an admissible tensor field under the change of  an
  admissible coordinate system  transform in the following way:
$$t^{a_{1},...,a_{p}}_{b_{1},...,b_{q}}=t^{\tilde{a}_{1},...,\tilde{a}_{p}}_{\tilde{b}_{1},...,\tilde{b}_{q}}A^{a_{1}}_{\tilde{a}_{1}}\cdots A_{b_{q}}^{\tilde{b}_{q}},$$
where $A^{a_{i}}_{\tilde{a}_{i}}=\frac{\partial
x^{a^{i}}}{\partial x^{\tilde{a}^{i}}}$.

The invariant character of the above statement is enclosed in the
equality
$$L_{\vec{\xi}}t^{a_{1},...,a_{p}}_{b_{1},...,b_{q}}=\partial_{n}t^{a_{1},...,a_{p}}_{b_{1},...,b_{q}},$$
where $L_{\vec{\xi}}$ is the Lie derivative along a vector field
$\vec{\xi}$.

We say that an admissible tensor field is integrable if there
exists an atlas of adapted coordinate maps such that the
components of this tensor in any of these coordinate maps are
constant. From Theorem \ref{Th1} immediately follows that the
necessary condition of the integrability of an admissible tensor
field  $t$ is vanishing of the derivatives  $\partial_{n}t$. We
call an admissible tensor structure $t$ quasi-integrable if in the
adapted coordinates it holds  $\partial_{n}t=0$. The form
$\omega=d\eta$ is an important example of an integrable admissible
structure. The following two theorems show the importance of the
just now given definitions.

\begin{theorem}\label{Th2}  The field of endomorphism  $\varphi$ is integrable if and only if $P(N_{\varphi})=0$ holds.
\end{theorem}

\textbf{Proof.} $\Rightarrow:$ The expression of the Nijenhuis
torsion
$$N_{\varphi}(\vec{X}, \vec{Y})=[\varphi \vec{X}, \varphi
\vec{Y}]+\varphi^{2}[\vec{X},\vec{Y}]-\varphi[\varphi\vec{X},\vec{Y}]-\varphi[\vec{X},\varphi\vec{Y}]$$
of the tensor $\varphi$ in adapted coordinates has the form:
\begin{align}\label{eq1} N^{e}_{ab}&=\varphi^{c}_{a}\vec{e}_{c}\varphi^{e}_{b}-
\varphi^{c}_{b}\vec{e}_{c}\varphi^{e}_{a}+\varphi^{e}_{c}\vec{e}_{b}\varphi^{c}_{a}-\varphi^{e}_{d}\vec{e}_{a}\varphi^{d}_{b},\\
\label{eq2}
N^{n}_{ab}&=2\varphi^{c}_{a}\varphi^{d}_{b}\omega_{dc}, \\
\label{eq3}
N^{e}_{na}&=-\varphi^{e}_{c}\partial_{n}\varphi^{c}_{a}, \\
\label{eq4} N^{n}_{na}&=0,\\
\label{eq5} N^{a}_{nn}&=0.\end{align} If the structure $\varphi$
is integrable, then from   \eqref{eq1}--\eqref{eq5} it follows
that
$$N_{\varphi}(\vec{e}_{a},\vec{e}_{b})=\varphi^{c}_{a}\varphi^{d}_{b}M^{n}_{cd}\partial_{n},\quad
N_{\varphi}(\partial_{n},\vec{e}_{a})=-(\partial_{n}\varphi^{c}_{b})\varphi^{a}_{c}\vec{e}_{a}.$$
The last two equalities imply  $P(N_{\varphi})=0$.

$\Leftarrow:$ Suppose that $P(N_{\varphi})=0$. Consider
sufficiently small neighborhood  $U$  of an arbitrary point of the
manifold $X$. Assume that  $U=U_{1} \times U_{2}$, $TU={\rm
span}(\partial_{a}) \oplus {\rm span}(\partial_{n})$. We set the
natural denotation $T(U_{1})={\rm span}(\partial_{a})$. We define
over the set   $U$ the isomorphism of bundles $\psi:D\rightarrow
T(U_{1})$ by the formula $\psi(\vec{e}_{a})=\partial_{a}$. This
endomorphism induces an almost complex structure on the manifold
$U_{1}$. This complex structure is integrable due to the equality
$P(N_{\varphi})=0$. Indeed, from \eqref{eq3} it follows that the
right hand side part of \eqref{eq1} coincides with the torsion of
the almost complex structure induced on the manifold $U_{1}$.
Choosing an appropriate coordinate system on $U_{1}$, and
consequently an appropriate adapted coordinate system on the
manifold  $X$, we get a coordinate map with respect to that the
components of the endomorphism  field $\varphi$ are constant.
$\Box$

\begin{theorem}\label{Th3} An almost contact metric structure is normal if and only
if the following conditions hold: $$P(N_{\varphi})=0,\quad
\omega(\varphi\vec{u},\varphi\vec{v})=\omega(\vec{u},\vec{v}).$$
\end{theorem}

\textbf{Proof.} Using the coordinate form
\eqref{eq1}--\eqref{eq5}, we see that the condition
$N_{\varphi}+2d\eta \otimes \vec{\xi}=0$ is equivalent to the
following system of equations: \begin{align}\nonumber
\varphi^{c}_{a}\vec{e}_{c}\varphi^{e}_{b}-\varphi^{d}_{b}\vec{e}_{c}\varphi^{e}_{a}+\varphi^{e}_{c}\vec{e}_{b}\varphi^{c}_{a}
-\varphi^{e}_{d}\vec{e}_{a}\varphi^{d}_{b}&=0,\\ \label{eq6}
-\varphi^{e}_{c}\partial_{n}\varphi^{c}_{a}&=0, \\
\nonumber
2\varphi^{c}_{a}\varphi^{d}_{b}\omega_{dc}&=2\omega_{ba}.\end{align}
This proves the theorem. $\Box$

The next statements shows the advisability of the notions like an
adapted coordinate system and an integrable tensor field.

\begin{theorem}\label{Th4} A contact metric structure is normal if and only if the field of endomorphisms  $\varphi$
is integrable. \end{theorem}

\textbf{Proof.} This statement follows from the fact that for a
contact metric structure the condition
$N_{\varphi}+2d\eta\otimes\vec{\xi}=0$ is equivalent to the
equality  $P(N_{\varphi})=0$, since the condition \eqref{eq6},
written in the coordinate-free form
$\omega(\varphi\vec{u},\varphi\vec{v})=\omega(\vec{u},\vec{v})$
holds automatically due to the definition of the contact metric
structure. $\Box$

Theorem \ref{Th4}  confirms the importance of the introducing of
the new types of the almost contact metric spaces. Namely, we call
an almost contact metric space {\it a Hermitian almost contact
metric space} if the condition $P(N_{\varphi})=0$ holds.

An intrinsic linear connection on a nonholonomic manifold $D$ is
defined in \cite{Lit10} as a map  $$\nabla:\Gamma D \times \Gamma
D \rightarrow \Gamma D $$ that satisfy the following conditions:
\begin{align*}
1)\quad &
\nabla_{f_1\vec{u}_1+f_2\vec{u}_2}=f_1\nabla_{\vec{u}_1}+f_2\nabla_{\vec{u}_2};\\
2)\quad &
\nabla_{\vec{u}}f\vec{v}=f\nabla_{\vec{u}}\vec{v}+(\vec{u}f)\vec{v},
\end{align*} where $\Gamma D$ is the module of admissible vector fields. The
Christoffel symbols are defined by the relation $$
\nabla_{\vec{e}_{a}}\vec{e}_{b}=\Gamma^{c}_{ab}\vec{e}_{c}. $$

The torsion $S$ of the intrinsic linear connection is defined by
the formula
$$
S(\vec{X},\vec{Y})=\nabla_{\vec{X}}\vec{Y}-\nabla_{\vec{Y}}\vec{X}-P[\vec{X},\vec{Y}].
$$ Thus with respect to an adapted coordinate system it holds  $$
S^{c}_{ab}=\Gamma^{c}_{ab}-\Gamma^{c}_{ba}. $$  In the same way as
a linear connection on a smooth manifold, an intrinsic connection
can be defined by giving a horizontal distribution over a total
space of some vector bundle. The role of such bundle plays the
distribution $D$.

In order to define a connection over the distribution $D$, it is
necessary first to introduce a structure of a smooth manifold on
$D$. This structure is defined in the following way. To each
adapted coordinate map  $K(x^\alpha)$ on the manifold  $X$ we put
in correspondence the coordinate map
$\tilde{K}(x^{\alpha},x^{n+\alpha})$ on the manifold  $D$, where
$x^{n+\alpha}$ are the coordinates of an admissible vector with
respect to the basis
$\vec{e}_{a}=\partial_{a}-\Gamma^{n}_{a}\partial_{n}$.

The notion of   {\it a connection over a distribution} introduced
in \cite{Lit12,Lit13}, was applied later to nonholonomic manifolds
with admissible Finsler metrics in \cite{Lit14,Lit15}. One says
that over a distribution  $D$ a connection is given if the
distribution $\tilde{D}=\pi^{-1}_{*}(D)$, where $\pi:D \rightarrow
X$ is the natural projection, can be decomposed into a direct some
of the form $$\tilde{D}=HD \oplus VD,$$ where $VD$ is the vertical
distribution on the total space $D$. Thus the assignment of a
connection over a distribution is equivalent to the assignment of
the object $G^{a}_{b}(X^{a},X^{n+a})$  such that
$$HD={\rm span}(\vec{\epsilon}_{a}),$$ where
$\vec{\epsilon}_{a}=\partial_{a}-\Gamma^{n}_{a}\partial_{n}-G^{b}_{a}\partial_{n+b}$.

It can be checked in the usual way that that the connection over
the distribution $D$  coincides with the linear connection in the
nonholonomic manifold $D$ if it holds
$$G^{a}_{b}(x^{a},x^{n+a})=\Gamma^{a}_{bc}(x^{a})x^{n+c}.$$ In \cite{Lit15} the notion of the prolonged connection was introduced. The prolonged
 connection can be obtained from an intrinsic connection by the equality
$$TD=\tilde{HD} \oplus VD,$$ where $HD \subset \tilde{HD}$. Essentially, the prolonged connection is a connection in a vector bundle.

An important example of a manifold with an admissible tensor
structure and a compatible with it intrinsic connection considered
V.V.~Wagner in \cite{Lit10}. In this paper, in a nonholonomic
manifold an intrinsic metric is introduced using an admissible
tensor field $g$ that satisfies the usual properties of the metric
tensor in a Riemannian space.

Similarly to the holonomic case, a metric on a nonholonomic
manifold defines there an intrinsic linear symmetric connection.
The corresponding Christoffel symbols can be derived from the
system of equations
$$
\nabla_{c}g_{ab}=\vec{e}_{c}g_{ab}-\Gamma^{d}_{ca}g_{db}-\Gamma^{d}_{cb}g_{ad}.
$$

Let  $\varphi$ be an admissible almost complex structure. We will
use the following statement.

\begin{theorem}\label{Th5}  Each nonholonomic manifold with an almost complex structure
 $\varphi$ and an intrinsic torsion-free linear connection
$\nabla$ admits an intrinsic linear connection $\tilde{\nabla}$
compatibel with the structure
  $\varphi$ and having the torsion $S$ such that
$$
S(\vec{u},\vec{v})=\frac{1}{4}PN_{\varphi}(\vec{u},\vec{v}),
$$
where $\vec{u},\vec{v}\in\Gamma(D)$.\end{theorem}

The proof of this theorem is based on the ideas of the proof of
Theorem 3.4 from \cite[Ch.~9]{Lit16}. One defines the following
connection $\tilde{\nabla}$: \begin{equation}\label{eq7}
\tilde{\nabla}_{\vec{u}}\vec{v}=\nabla_{\vec{u}}\vec{v}-Q(\vec{u},\vec{v}),\end{equation}
where $Q$ can be constructed in a special way using  $\nabla$ and
$\varphi$ \cite[P.~137]{Lit16}.

\section{Interior characteristics of almost contact metric spaces}

Now we introduce the notion of an almost contact metric structure
in a new sense. Namely, we will say that a manifold of almost
contact metric structure in the new sense is given if on a
manifold $X$ with a given contact form  $\eta$, in addition a pair
of admissible tensor structures   $(\varphi, g)$ such that
$$\varphi^{2}\vec{u}=-\vec{u},\quad
g(\varphi\vec{u},\varphi\vec{v})=g(\vec{u},\vec{v})$$ is given.

\begin{theorem}\label{Th6}
 The notion of a manifold of almost contact
metric structure in the new sense is equivalent to the notion of a
manifold of almost contact metric structure in the old sense.
\end{theorem}

{\bf Proof.} Suppose that on the manifold  $X$ an  almost contact
metric structure in the new sense is given. We introduce on the
manifold  $X$ the Riemannian metric $\tilde{g}$ by the equality
$$\tilde{g}(\vec{u},\vec{v})=g(\vec{u},\vec{v}),$$ where
$\vec{u},\vec{v}\in \Gamma D$, $\tilde{g}(\vec{u},\vec{\xi})=0$
and $\tilde{g}(\vec{\xi},\vec{\xi})=1$. The necessary conditions
on field of endomorphisms $\varphi$ and the 1-form $\tilde{g}$ can
be checked directly.

We say that a Sasakian manifold in the new sense is given if on
the manifold $X$ with a given contact metric structure, in
addition the equality  $P(N)=0$ holds. Theorems \ref{Th3} and
\ref{Th6} imply the following statement.

\begin{theorem}\label{Th7} The notion of a Sasakian manifold in
the new sense is equivalent to the notion of a Sasakian manifold
in the old sense.\end{theorem}

In this section we use the following notation. As above,
admissible almost complex structure and Riemannian metric will be
denoted by  $\varphi$ and $g$, respectively; the symbol $\nabla$
will denote the intrinsic metric connection, and the symbols
$\tilde{g}$ and $\tilde{\nabla}$  will denote the metric tensor in
the ambient space and its Levi-Civita connection, respectively.

\begin{theorem}\label{Th8}
A contact metric structure is normal if and only if the structure
$\varphi$ is quasi-integrable and it holds  $\nabla\varphi=0$,
where $\nabla$ is an intrinsic metric connection. \end{theorem}

{\bf Proof.} $\Rightarrow:$ From Theorem~\ref{Th3} it follows that
 $P(N_{\varphi})=0$, and consequently the structure  $\varphi$ is quasi-integrable.
 Moreover, constructing, using the metric connection, the
 following connection satisfying
$$S(\vec{u},\vec{v})=\frac{1}{4}PN_{\varphi}(\vec{u},\vec{v}),\quad
\vec{u},\vec{v}\in\Gamma(D),$$ we get the second statement.

$\Leftarrow:$ From \eqref{eq7} it follows that
$Q(\vec{u},\vec{v})=0$ and $PN_{\varphi}(\vec{u},\vec{v})=0$. This
and the quasi-integrability of the structure  $\varphi$ imlies the
theorem. $\Box$

Note that the equality $\nabla\varphi=0$ is not true if the
connection
 $\nabla$ and the field of endomorphisms  $\varphi$  are considered as the structures defined on the whole manifold,
 see e.g. \cite{Lit7}.

Next suppose that $\nabla^{1}$ is the extended connection
constructed from the intrinsic connection in the following way:
$$\tilde{HD}=HD \oplus {\rm span}(\partial_{n}),$$ here
$\partial_{n}$ is a vector field on the manifold  $D$. The
extended connection allows to formulate the next characteristic
feature of the integrability of an almost complex structure
$\varphi$.

\begin{theorem}\label{Th9}
An almost complex structure  $\varphi$ is integrable if and only
if the equality $\nabla^{1}\varphi=0$ holds.\end{theorem}

Finally we formulate a statement concerning  $K$-contact
manifolds.

\begin{theorem}\label{Th10}
An almost contact metric structure is a $K$-contact structure if
and only if  the metric $g$ is quasi integrable.
\end{theorem}

The theorem follows from the following equivalences: $$
L_{\vec{\xi}}\tilde{g}=0 \Leftrightarrow L_{\vec{\xi}}g=0
\Leftrightarrow \partial_{n}g=0. $$

\bibliographystyle{unsrt}

\vskip2cm

Saratov State University,\\
Chair of Geometry\\
E-mail: sgalaev(at)mail.ru

\end{document}